\documentclass[twocolumn]{svjour3}
\usepackage{verbatim}  
\usepackage{graphicx}
\usepackage{amsfonts,amsmath,tabularx}

\newcommand{\bm}[1]{\mbox{\boldmath $#1$}}

\begin{document}

\title{Optimal design of fibre reinforced membrane structures}

\author{Anders~Klarbring \and Bo Torstenfelt \and Peter Hansbo \and Mats G. Larson}
\institute{A. Klarbring \at Division of Solid Mechanics, Link\"oping University, Sweden \email{anders.klarbring@liu.se}
\and B. Torstenfelt \at Division of Solid Mechanics, Link\"oping University, Sweden
\and P. Hansbo \at Department of Mechanical Engineering, J\"onk\"oping University, Sweden
\and M. G. Larson \at Department of Mathematics and Mathematical Statistics, Ume{\aa} University, Sweden}

\date{Received: date / Accepted: date}    

\maketitle
\begin{abstract}
A design problem of finding an optimally stiff membrane structure by selecting one--dimensio\-nal fiber reinforcements is formulated and solved. The membrane model is derived in a novel manner from a particular three-dimensional linear elastic orthotropic model by appropriate assumptions. The design problem is given in the form of two minimization statements, reminiscent of a Nash game. After finite element discretization, the separate treatment of each of the two minimization statements follows from classical results and methods of structural optimization: the stiffest orientation of reinforcing fibers coincides with principal stresses and the separate selection of density of fibers is a convex problem that can be solved by optimality criteria iterations. Numerical solutions are shown for two particular configurations. The first for a statically determined structure and the second for a statically undetermined one. The latter shows related but non-unique solutions.
\end{abstract}
\section{Introduction}
A finite element membrane shell model was recently derived by Hansbo and Larson \cite{HL} using tangential
differential calculus, meaning that the problem is set in a Cartesian three dimensional space as opposed to a
parametric plane, thereby generalizing the classical flat facet element shell model to higher order
elements. The present study further extends this membrane model by allowing for non-isotropic materials.
In particular, one--dimensional fibers are added to a base material, modeling, e.g., the reinforcements seen
in modern racing boat sails. The plane stress property, as well as the membrane property of complete out-of-plane
shear flexibility, is shown to be exact consequences of certain material parameter selections for a
three-dimensional transversely isotropic base material. This together with a displacement assumption
reduces the three dimensional model to the surface model. Based on this finite element
model we formulate a design problem where we seek to find the best fiber reinforcements of the
membrane, meaning that we find the stiffest structure by both rotation and sizing of the fibers. The formulation is reminiscent of a Nash game \cite{Aubin,HTK}, consisting of two minimization statements. However, the two players of the game have the same objective, i.e., stiffness, as opposed to standard game theory. Nevertheless, since the two minimization statements relates to rotation and sizing of the fibers, respectively, such a formulation ties directly to the sequential iterative treatment suggested for similar problems previously \cite{BS}. The optimal rotation is found by identifying the material as a so-called low shear material, implying that the optimal orthotropic principal directions coincides with the principal stress directions  \cite{Pedersen1,Pedersen2}, while the optimal thickness distribution is found by a classical optimality criteria iteration formula.

\section{The model}
We consider a material that is a mixture of a transversely isotropic linear elastic base material and $n$ reinforcing fibre materials. The transversely isotropic material has material constants that satisfy the plane stress assumption as well as the membrane behaviour of having complete flexibility when sheared perpendicularly to the membrane surface.

\subsection{Geometry}
The geometry of the membrane is defined by an orientable smooth surface $\Sigma$ with outward normal $\bm{n}$. For any point $\bm{x}\in \mathbb{R}^3$ we denote the signed distance function relative to $\Sigma$ by $\zeta(\bm{x})$. The membrane with thickness $t$ then occupies
\[
\Omega_t=\{\bm{x}\in  \mathbb{R}^3 : |\zeta(\bm{x})|<t/2\}.
\]
Note that $\nabla \zeta(\bm{x})=\bm{n}$ for $\bm{x}\in\Sigma$. For a sufficiently small $t$, the orthogonal projection point $\bm{p}(\bm{x})\in \Sigma$ of $\bm{x}\in \Omega_t$ is unique and given by
\[
\bm{p}(\bm{x})=\bm{x}-\zeta(\bm{x})\nabla \zeta(\bm{x}).
\]
Moreover, for $\bm{x}\in \Omega_t$, the linear projection operator of vectors onto the tangent plane of $\Sigma$ at $\bm{p}(\bm{x})$ is
\[
\bm{P}_\Sigma=\bm{I}-\bm{n}\otimes\bm{n},
\]
where $\bm{I}$ is the identity tensor and $\otimes$ denotes exterior product. In the sequel we will also need the projection operator onto the one-dimensional subspace spanned by $\bm{n}$, i.e.,
\[
\bm{N}_\Sigma=\bm{n}\otimes\bm{n}.
\]
Note that $\bm{P}_\Sigma\bm{N}_\Sigma=\bm{N}_\Sigma\bm{P}_\Sigma=\bm{0}$. The directions of the reinforcing fibers are given by vector fields $\bm{s}_i$, $i=1,\ldots, n$, such that $\bm{s_i}\cdot\bm{n}=0$. Projections onto these directions are then defined by
\[
\bm{S}_i=\bm{s}_i\otimes\bm{s}_i.
\]
Clearly $\bm{P}_\Sigma\bm{S}_i=\bm{S}_i\bm{P}_\Sigma=\bm{S}_i$.

\subsection{The material}
The base material is transversely isotropic with respect to an axis defined by $\bm{n}$. Such a material can be described by an elasticity tensor expressed in terms of five material constants $\delta_1$, $\delta_1$, $\delta_1$, $\gamma$ and $\mu$ according to \cite{LC,NPG}, so that the fourth order tensor of elastic moduli of the base material $\mathbb{E}^{\mbox{\scriptsize base}}$ can be written
\begin{align}\nonumber
\mathbb{E}^{\mbox{\scriptsize base}}= {}& \delta_1\bm{N}_\Sigma\otimes\bm{N}_\Sigma
+\delta_2(\bm{N}_\Sigma\otimes\bm{P}_\Sigma+\bm{P}_\Sigma\otimes\bm{N}_\Sigma)
\\ \nonumber {}&
+{\mu}(\bm{P}_\Sigma\underline{\otimes}\bm{P}_\Sigma+\bm{P}_\Sigma\overline{\otimes}\bm{P}_\Sigma)
+\delta_3\bm{P}_\Sigma\otimes\bm{P}_\Sigma+\\ \nonumber {}&
+\frac{\gamma}{2}(\bm{N}_\Sigma\underline{\otimes}\bm{P}_\Sigma+\bm{N}_\Sigma\overline{\otimes}\bm{P}_\Sigma \\
 &{} \qquad+\bm{P}_\Sigma\underline{\otimes}\bm{N}_\Sigma+\bm{N}_\Sigma\overline{\otimes}\bm{P}_\Sigma) .
\label{eq:E-trans}\end{align}
Here dyadic products of second order tensors are defined by their action on a third second order tensor, i.e.,
\begin{gather*}
(\bm{A}\otimes\bm{B})\bm{C}=(\bm{B}:\bm{C})\bm{A},\quad (\bm{A}\underline{\otimes}\bm{B})\bm{C}=\bm{A}\bm{C}\bm{B}^T,\\
(\bm{A}\overline{\otimes}\bm{B})\bm{C}=\bm{A}\bm{C}^T\bm{B}^T,
\end{gather*}
where a double dot indicates inner product of second order tensors.

The reinforcing fibers have elasticity tensors of the form
\begin{equation}\label{eq:E-fiber}
\mathbb{E}^{\mbox{\scriptsize fiber}}_i=\alpha_i\bm{S}_i\otimes\bm{S}_i, \quad 1=1,\ldots, n,
\end{equation}
where $\alpha_i$ are Young type elasticity coefficients.

The constitutive law of the membrane material is now taken as being composed of a constrained mixture of base material and reinforcing material. The amount of each material is defined by fractions $t_b$ and $t_i$, $i=1,\ldots, n$, of the membrane thickness $t$ such that the total constitutive tensor is given as
\[
\mathbb{E}=\frac{t_b}{t}\mathbb{E}^{\mbox{\scriptsize base}}+\sum_{i=1}^n\frac{t_i}{t}\mathbb{E}^{\mbox{\scriptsize fiber}}_i, \quad t=t_b+\sum_{i=1}^nt_i,
\]
and the linear constitutive law is
\begin{equation}\label{eq:constitutive}
  \bm{\sigma}=\mathbb{E}[\bm{\varepsilon}],
\end{equation}
where $\bm{\sigma}$ and $\bm{\varepsilon}$ are the stress and strain tensors, respectively.

\subsection{Membrane stress assumptions}
We define a membrane material by the requirements that it is always in a state of plane stress and no shear stress perpendicular to the membrane surface exists, i.e.,
\begin{equation}\label{eq:membrane}
  \bm{N}_\Sigma\bm{\sigma}\bm{N}_\Sigma=\bm{0},\quad \bm{P}_\Sigma\bm{\sigma}\bm{N}_\Sigma=\bm{N}_\Sigma\bm{\sigma}\bm{P}_\Sigma=\bm{0}.
\end{equation}
The zero bending stiffness behaviour of membranes will be a result of a kinematic assumption introduce subsequently.
Inserting (\ref{eq:constitutive}) into (\ref{eq:membrane}) gives
\begin{align}\nonumber
  \bm{N}_\Sigma\bm{\sigma}\bm{N}_\Sigma= {}& \frac{t_b}{t}[\delta_1\bm{N}_\Sigma(\bm{N}_\Sigma:\bm{\varepsilon})+
  \delta_2\bm{N}_\Sigma(\bm{P}_\Sigma:\bm{\varepsilon})] \\
  = {}& \bm{0},\label{eq:normal_stress}
\end{align}
\begin{equation}\label{eq:shear_stress}
  \bm{P}_\Sigma\bm{\sigma}\bm{N}_\Sigma=\frac{t_b}{t}[\gamma\bm{P}_\Sigma\bm{\varepsilon}\bm{N}_\Sigma]=\bm{0}.
\end{equation}
Thus, we conclude that the constitutive constant $\gamma$ needs to be zero and that the strain perpendicular to the membrane is controlled by the in-plane strain as
\begin{equation}\label{eq:plane_stress}
  \bm{N}_\Sigma:\bm{\varepsilon}=-\frac{\delta_2}{\delta_1}\bm{P}_\Sigma:\bm{\varepsilon}.
\end{equation}
Moreover, the in-plane stress can be calculated from (\ref{eq:constitutive}) as follows:
\begin{align*}
  \bm{P}_\Sigma\bm{\sigma}\bm{P}_\Sigma=\frac{t_b}{t}[ & \delta_2\bm{P}_\Sigma(\bm{N}_\Sigma:\bm{\varepsilon})
  +\delta_3\bm{P}_\Sigma(\bm{P}_\Sigma:\bm{\varepsilon}) \\
 & +2\mu\bm{P}_\Sigma\bm{\varepsilon}\bm{P}_\Sigma]+\sum_{i=1}^n\frac{t_i}{t}\alpha_i\bm{S}_i(\bm{S}_i:\bm{\varepsilon}),
\end{align*}
and when using (\ref{eq:plane_stress}) we get
\begin{align}\nonumber
  \bm{P}_\Sigma\bm{\sigma}\bm{P}_\Sigma= {}&
  \frac{t_b}{t}[\delta\bm{P}_\Sigma(\bm{P}_\Sigma:\bm{\varepsilon})+2\mu\bm{P}_\Sigma\bm{\varepsilon}\bm{P}_\Sigma]\\ &{}+\sum_{i=1}^n\frac{t_i}{t}\alpha_i\bm{S}_i(\bm{S}_i:\bm{\varepsilon}),
\label{eq:plane_stress2}\end{align}
where
\[
\delta=\delta_3-\frac{\delta_2^2}{\delta_1}.
\]
The elasticity coefficient $\mu$ equals the in-plane shear modulus, while $\delta$ is a plane stress Lam\'{e} coefficient. The two elasticity moduli $\delta$ and $\mu$ can be expressed in terms of in-plane Young and Poisson moduli $E$ and $\nu$ as
\[
\delta=\frac{\nu E}{1-\nu^2},\quad \mu = \frac{E}{2(1+\nu)}.
\]

The volumetric specific strain energy can, due to (\ref{eq:membrane}) be written as
\[
W_s=\frac{1}{2}\bm{\sigma}:\bm{\varepsilon}=
\frac{1}{2}(\bm{P}_\Sigma\bm{\sigma}\bm{P}_\Sigma):\bm{\varepsilon}.
\]
Inserting (\ref{eq:plane_stress2}) we get
\[
W_s=\frac{1}{2}\left(\mathbb{E}^{\mbox{\scriptsize memb}}[\bm{\varepsilon}]\right):\bm{\varepsilon},
\]
where the membrane elasticity tensor is defined by
\begin{align*}
\mathbb{E}^{\mbox{\scriptsize memb}}= {}& \frac{t_b}{t}[\delta\bm{P}_\Sigma\otimes\bm{P}_\Sigma+\mu(\bm{P}_\Sigma\underline{\otimes}\bm{P}_\Sigma+\bm{P}_\Sigma\overline{\otimes}\bm{P}_\Sigma)] \\ {}& +
\sum_{i=1}^n\frac{t_i}{t}\alpha_i\bm{S}_i\otimes\bm{S}_i.
\end{align*}

\subsection{Potential energy}
The strain is derived as usual as the symmetrized gradient of the displacement vector $\bm{u}$:
\[
\bm{\varepsilon}=\bm{\varepsilon}(\bm{u})=\frac{1}{2}(\nabla\otimes\bm{u}+(\nabla\otimes\bm{u})^\text{T}).
\]
Therefore, we can regard the volume specific strain energy as a function of the displacement field, i.e., $W_s=W_s(\bm{u})$.

We now introduce the basic kinematic assumption that all material points in $\Omega_t$ that lie along a normal to the surface $\Sigma$ have the same displacement vector, i.e.,
\[
\bm{u}(\bm{x})=\bm{u}(\bm{p}(\bm{x})),\quad \bm{x}\in\Omega_t.
\]
This kinematics imply that bending of the membrane is essentially eliminated and no bending stiffness, despite the finite thickness, is present.

The total strain energy, which is the volume integral of $W_s$ can then be written:
\[
W=\int_{-t/2}^{t/2}\int_\Sigma W_s(\bm{u}(\bm{p}(\bm{x}))\;d\Sigma_\zeta d\zeta,
\]
where $d\Sigma_\zeta$ is an area element for a surface parallel to $\Sigma$ at the distance $\zeta$, which reads
\[
d\Sigma_\zeta=d\Sigma (1+\zeta H+\zeta^2 K),
\]
where $d\Sigma$ is the area element of $\Sigma$, and $H$ and $K$ are the mean curvature and Gaussian curvature, respectively. For a membrane that is thin compared to its curvature we can use the approximation
\[
d\Sigma_\zeta\approx d\Sigma.
\]
The total potential energy is now taken as
\[
\Pi=\frac{t}{2}\int_\Sigma W_s(\bm{u}(\bm{x}))\;d\Sigma-\langle \bm{F},\bm{u}\rangle_\Sigma,
\]
where the force $\bm{F}$ is a member of the dual space of displacement fields on $\Sigma$ and $\langle \cdot,\cdot\rangle_\Sigma$ is a duality paring.

\section{Equilibrium}
We define the membrane forces (per unit length) as
\[
\bm{M}=t\bm{P}_\Sigma\bm{\sigma}\bm{P}_\Sigma.
\]
Stationarity of the potential energy gives the following principle of virtual work:
\begin{equation}\label{eq:virtual_work}
  \int_\Sigma\bm{M}:\bm{\varepsilon}(\bm{v})\;d\Sigma=\langle \bm{F},\bm{v}\rangle_\Sigma,
\end{equation}
for all kinematically admissible fields $\bm{v}$. Such fields will generally be restricted in the tangential direction on a subset of $\partial\Sigma$. We will assume that loading on the membrane can be written as
\[
\langle \bm{F},\bm{v}\rangle_\Sigma=\int_\Sigma\bm{f}\cdot\bm{v}\;d\Sigma+\int_S\bm{p}\cdot\bm{v}\;dS,
\]
where $\bm{f}$ is a force per area over $\Sigma$, and $\bm{p}$ is a force per unit length over the part $S$ of $\partial\Sigma$ where the displacement is not prescribed. Using now Lemma 2.1 of Gurtin and Murdoch \cite{GM}, i.e., an integral theorem for surfaces, we obtain the equilibrium equations
\begin{equation}\label{eq:equilibrium1}
  -\mbox{div}_\Sigma\:\bm{M}=\bm{f},
\end{equation}
\begin{equation}\label{eq:equilibrium2}
  \bm{M}\bm{\nu}=\bm{p},
\end{equation}
where $\mbox{div}_\Sigma$ is the surface divergence, and $\bm{\nu}$ is a unit vector of $\partial\Sigma$, tangential to $\Sigma$. Since $\bm{M}\bm{\nu}$ will also be a vector tangent to $\Sigma$ we conclude that $\bm{p}$ can have no component perpendicular to the surface.

\section{Design problem}
From now on we will consider the special case of an orthotropic material consisting of two orthogonal families of fibers, consisting of the same material, i.e., $\alpha_1=\alpha_2=\alpha$. From now on we use the notation $\bm{s}=\bm{s}_1$ and $\bm{s}^\bot=\bm{s}_2$.

The orientation of the fibers in the tangent plane of the membrane, i.e., $\bm{s}$ and $\bm{s}^\bot$, can be defined by an angle $\theta$ belonging to
\[
\Theta=\{\theta|\;0\leq\theta\leq 2\pi\}.
\]
This angle will be a design variable in the optimal design problem. Other such design variables are $t_1$ and $t_2$, i.e., the fiber contents in the two orthogonal directions. The field $\bm{t}=(t_1,t_2)$ belongs to the set
\begin{align*}
T=\Bigl\{ & \bm{t}=(t_1,t_2)|\;\underline{t}_\alpha\leq t_\alpha\leq\overline{t}_\alpha, \alpha=1,2,\;\\
& \int_{\Sigma}(t_1+t_2)d\Sigma\leq V\Bigr\},
\end{align*}
where $\underline{t}_\alpha$ and $\overline{t}_\alpha$ are non-negative upper and lower bounds and $V$ is a limit for the total amount of material that can be used for the fibers.

The potential energy is seen as a function
\[
\Pi:V\times T\times \Theta \rightarrow \mathbb{R},
\]
where $V$ is the set of kinematically admissible displacements. Minimizing $\Pi$ with respect to the first argument gives the equilibrium displacement as a function of the design variables, i.e., $\bm{u}=\bm{u}(\bm{t},\theta)$. As a measure of stiffness we use the so called compliance
\begin{align*}
C(\bm{t},\theta):= {}& \frac{1}{2}\langle \bm{F},\bm{u}(\bm{t},\theta)\rangle_\Sigma=-\Pi(\bm{u}(\bm{t},\theta), \bm{t},\theta)\\
= {}&-\min_{v\in V}\Pi(\bm{v},\bm{t},\theta).
\end{align*}
Our design goal is to find a design that minimizes the compliance. We choose to split this into two parts as follows: find $\bm{t}^*\in T$ and $\theta^* \in \Theta$ such that
\[
(\mathbb{P})\quad\left\{
\begin{array}{l}
  C(\bm{t}^*,\theta^*)=\min_{\bm{t}\in T}C(\bm{t},\theta^*)\\[5pt]
  C(\bm{t}^*,\theta^*)=\min_{\theta \in \Theta}C(\bm{t}^*,\theta).
\end{array}\right.
\]
This is reminiscent of a Nash equilibrium problem \cite{Aubin,HTK}, but the two players have the same objective in contrast to natural Nash games where objectives are opposing, which in the present case would happen if minimization is changed for maximization in one of the two lines of $(\mathbb{P})$.

The second sub-problem of $(\mathbb{P})$, i.e., finding an optimal orientation for an orthotropic material, has been extensively discussed by Pedersen~\cite{Pedersen1,Pedersen2} and Hammer~\cite{Hammer}, see also Bendsoe and Sigmund~\cite{BS} for further discussion and related references. It turns out that the problem can be solved locally, i.e., the orientation of the material is determined by the local stress state only, and in particular the orientation of principal stresses and strains. Due to the plane stress assumption there are only two possibly non-zero principal components of the stress tensor $\bm{\sigma}$, denoted $\sigma_I$ and $\sigma_{II}$, such that $|\sigma_I|\geq |\sigma_{II}|$. The corresponding principal directions (eigenvectors) are tangent to the membrane plane. Obviously, these facts also holds for the principal components of $\bm{M}$, i.e., $M_I$ and $M_{II}$, such that $|M_I|\geq |M_{II}|$. For a so-called low shear orthotropic material, the solution $\theta^*$ of the second sub-problem of $(\mathbb{P})$ represents an orientation where the orthotropic principal directions coincide with the principal stress or membrane force directions, which are also the principal strain directions. Moreover, the orthotropic principal direction having that highest stiffness should be in the direction corresponding to $\sigma_I$ and $M_I$. In the Appendix we show that the particular orthotropic material defined above, having by two families of fibers in orthogonal directions $\bm{s}$ and $\bm{s}^\bot$, is indeed a low shear material and, therefore, the optimal directions of $\bm{s}$ and $\bm{s}^\bot$ are in the directions of principal stress. Moreover, if $t_1 > t_2$ then $\bm{s}$ is in the direction of $\sigma_I$.

The first sub-problem of $(\mathbb{P})$ is a classical stiffness optimization problem, albeit having two design fields, one for each fiber orientation. This is a convex problem and can be solved by satisfying the optimality conditions. The surface elasticity tensor $\mathbb{S}^{\mbox{\scriptsize memb}}=t\,\mathbb{E}^{\mbox{\scriptsize memb}}$ is regarded as a function of the design, i.e., $\mathbb{S}^{\mbox{\scriptsize memb}}=\mathbb{S}^{\mbox{\scriptsize memb}}(\bm{t},\theta)$. The optimality conditions of the first sub-problem of $(\mathbb{P})$ become \cite{BS,CK}:
\begin{equation}\label{eq:opt1}
  \frac{\partial \mathbb{S}^{\mbox{\scriptsize memb}}}{\partial t_\alpha}[\bm{\varepsilon}(\bm{u})]:\bm{\varepsilon}(\bm{u})=\Lambda+\lambda_\alpha^+-\lambda_\alpha^-,\;\;\alpha=1,2,
\end{equation}
\begin{equation}\label{eq:opt2}
  \Lambda\geq 0,\quad \Lambda\left(\int_{\Sigma}(t_1+t_2)d\Sigma-V\right)=0,
\end{equation}
\begin{gather}\nonumber
  \lambda_\alpha^+\geq 0,\;\lambda_\alpha^-\geq 0,\;\lambda_\alpha^-(\underline{t}_\alpha-t_\alpha)=0,\\
  \;\lambda_\alpha^+({t}_\alpha-\overline{t}_\alpha)=0,\;\;\alpha=1,2,
\label{eq:opt3}\end{gather}
where $\Lambda$, $\lambda_\alpha^+$ and $\lambda_\alpha^-$ are Lagrangian multipliers, $\bm{t}\in T$ and $\bm{u}=\bm{u}(\bm{t},\theta)$ is the displacement solution, i.e., the minimum field with respect to $\bm{v}$ of $\Pi(\bm{v},\bm{t},\theta)$.

Note that
\[
\frac{\partial \mathbb{S}^{\mbox{\scriptsize memb}}}{\partial t_\alpha}={\alpha}\bm{S}_\alpha\otimes\bm{S}_\alpha.
\]

\section{Discretization and algorithm}
For the numerical treatment of $(\mathbb{P})$ we need to introduce a discrete approximation. The discretization of the state problem, i.e., the problem of finding the minimum displacement $\bm{u}\in V$ of the potential energy $\Pi$ for a given design $\theta\in\Theta$ and $\bm{t}\in T$, follows Hansbo and Larson \cite{HL}. This implies introducing a triangulation of $\Sigma$ resulting in a discrete surface, with corresponding discrete normal vector field and projections. The displacement field is approximated using the same triangulation but is possibly of different order.

In addition to the approximation of the state problem we also need to approximate the design fields $\bm{t}\in T$ and $\theta\in \Theta$. This is achieved by using point values: these are denoted $\bm{t}_i=(t_{1i},t_{2i})$ and $\theta_i$ for point $i$. In particularly, we use superconvergence points of the finite elements \cite{Barlow}. Such a discretization means that (\ref{eq:opt1}) and (\ref{eq:opt3}) are imposed at these evaluation points and the integral in (\ref{eq:opt2}) is replaced by a sum.

Let
\[
A^k_{\alpha i}=\left( \frac{\partial \mathbb{S}^{\mbox{\scriptsize memb}}}{\partial t_\alpha}[\bm{\varepsilon}(\bm{u^k})]:\bm{\varepsilon}(\bm{u^k})\right)_i,
\]
be the left hand side of (\ref{eq:opt1}) evaluated at point $i$ and for a displacement field $\bm{u^k}$. Also, let $B^k_{\alpha i}=(\Lambda^k)^{-1}A^k_{\alpha i}$ where $\Lambda^k$ is a currant iterate of the Lagrangian multiplier $\Lambda$. For a given displacement iterate $\bm{u^k}$ and rotation $\theta^k$ the following fixed point iteration formula is suggested by the optimality conditions (\ref{eq:opt1}) through (\ref{eq:opt3}):
\begin{equation}\label{eq:OC}
  t^{k+1}_{\alpha i}=\left\{\begin{array}{ll}
                   \underline{t}_{\alpha i} & \mbox{if } t^k_{\alpha i} (B^k_{\alpha i})^\eta\leq \underline{t}_{\alpha i} \\[2mm]
                   \overline{t}_{\alpha i} & \mbox{if } t^k_{\alpha i} (B^k_{\alpha i})^\eta\geq \overline{t}_{\alpha i} \\[2mm]
                   t^k_{\alpha i} (B^k_{\alpha i})^\eta & \mbox{otherwise},
                 \end{array}\right.
\end{equation}
where $\underline{t}_{\alpha i}$ and $\overline{t}_{\alpha i}$ are point values of the upper and lower bounds and $0<\eta\leq 1$  is a damping coefficient.

The following algorithmic steps, the convergence of which gives satisfaction of a discrete version of the optimality conditions of $(\mathbb{P})$, is now suggested:
\begin{enumerate}
  \item For a given design $\theta^k$ and $\bm{t}^k$, solve the state problem, i.e., find the minimum displacement field of $\Pi(\bm{v},\bm{t}^k ,\theta^k)$ so as to obtain the currant displacement iterate $\bm{u}^k$.

  \item Obtain new fiber thickness distributions by the optimality criteria formula (\ref{eq:OC}) where
  \begin{itemize}
    \item $\Lambda_k$ is determined such that
    \[
    \sum_i(t_{1i}^{k+1}+t_{2i}^{k+1})\:d\Sigma\leq V.
    \]
    A local iteration is needed for this.
  \end{itemize}

  \item For each integration point, calculate principal stresses (and/or principal membrane forces). Take $\bm{s}$ to correspond to the main material direction, i.e., to $t_{1i}$, such that $t_{1i}\geq t_{2i}$, and chose $\theta^{k+1}$  so that this $\bm{s}$ aligns with the main principal stress direction.

  \item Let $k=k+1$ and return to the first step.
\end{enumerate}
Steps 1 and 2 can be iterated several times before continuing with calculation of fiber directions in Step 3. In fact, in the examples the fixed point iteration (\ref{eq:OC}), for newly calculated displacement $\bm{u}^k$, is repeated until convergence before continuing with the fiber directions in Step 3.

Note that step 3 assumes distinct principal stresses. Numerically coalescence of such stresses occur with close to zero probability but may show up as non-convergence issues. For statically determined structures, i.e., when $\bm{M}$ is uniquely determined by (\ref{eq:equilibrium1}) and (\ref{eq:equilibrium2}), this may be of particular concern. For such cases that have distinct principal stresses, step 3 above needs to be performed only ones since these principal stress are independent of $\bm{t}$. Such problems essentially becomes convex since the first part of $(\mathbb{P})$ is a convex problem. The first problem of the Section \ref{sec:examples} is statically determinate but has not everywhere distinct principal stresses.

\section{Examples}\label{sec:examples}
\subsection{Oblate spheroid}
An oblate spheroid, where $\Sigma$ is defined by
\[
x^2+y^2+(2z)^2=1,
\]
was solved by different finite elements and triangulations in Hansbo and Larson \cite{HL}. Here we treat the same geometry but use an internal pressure $p$ as loading. We seek for optimal fiber distribution as described in previous sections. The data are $E=1$, $\nu=0.3$, $t_b=0.005$, $p=10$, $V=0.01$, $\overline{t}_1=\overline{t}_2=0.004$, $\underline{t}_1=\underline{t}_2=0$ and $\alpha=1$. The initial fiber thickness is uniform and chosen so as to satisfy the volume constraint as an equality. We use 3072 bilinear 4-node fully integrated isoparametric elements, implying one superconvergent point per element and, thus, three design variables per element. Symmetry is utilized and only half of the spheroid is modeled. The problem converged in 36 optimality criteria updates and 7 updates of the fiber orientations. As convergence criteria an objective value change below 0.001 \% and a change of $\theta$ such that $\cos\theta > 0.999$ are used. Note that the problem is statically determinate, but at the poles of the spheroid symmetry implies that the principal stresses coincide for an exact solution. This is the reason for the need of several updates of fiber orientations before convergence, despite the problem being statically determinate.
\begin{figure}[htp]
\centering
\includegraphics[width=7cm]{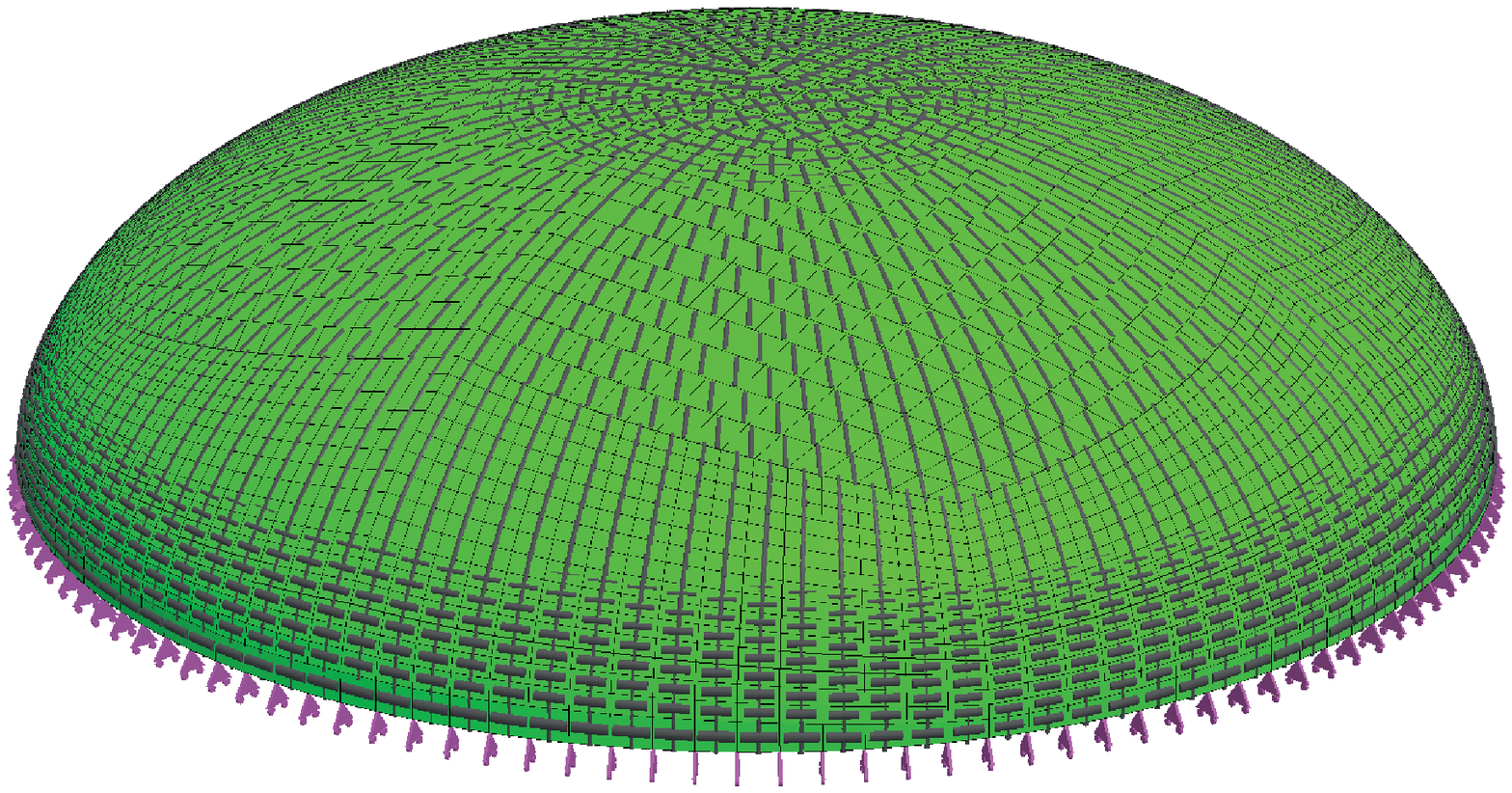}
\includegraphics[width=7cm]{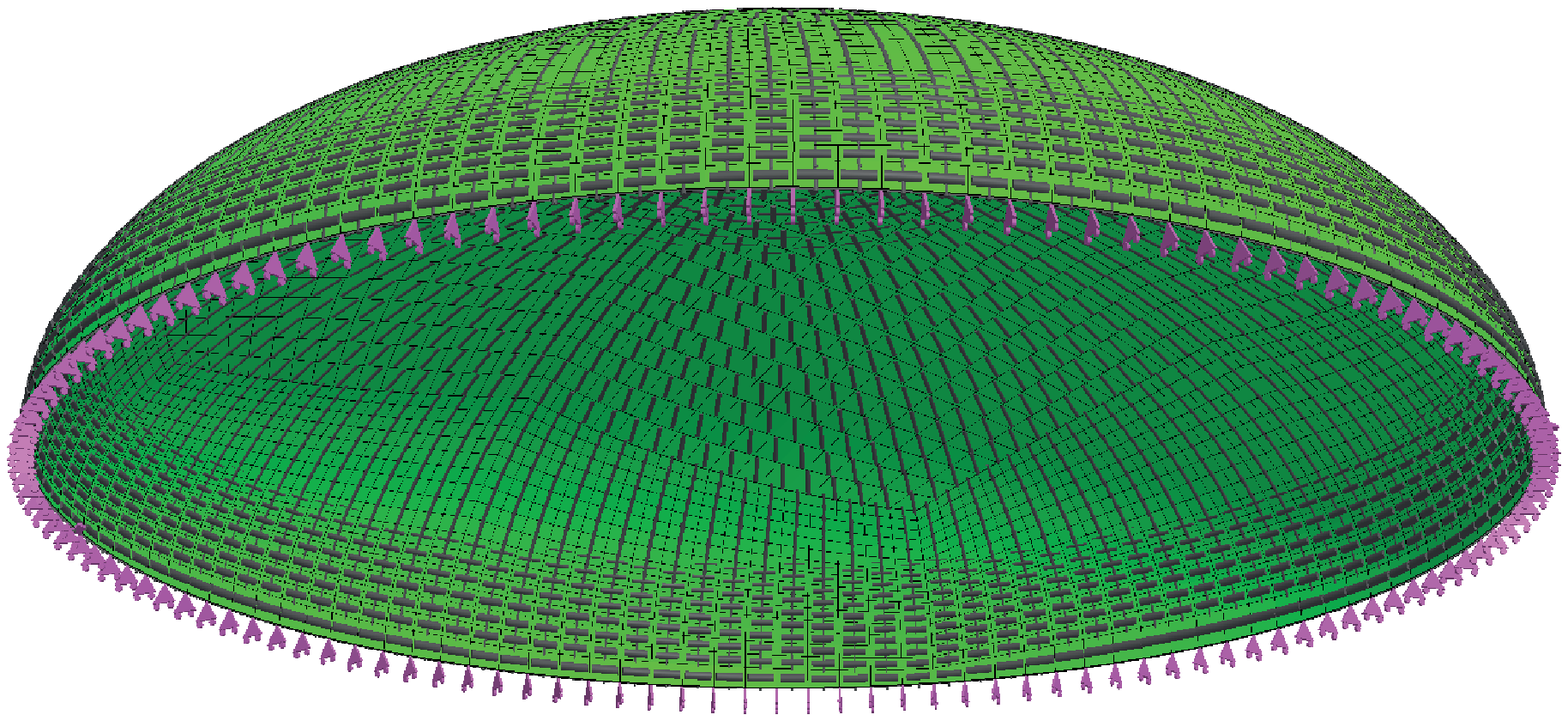}
\caption{Optimal fiber distribution of an oblate spheroid, loaded by internal pressure.}
\label{fig:example1}
\end{figure}

What concerns the general features of the solution one finds, on examination of Figure \ref{fig:example1}, that close to the equator both fiber families are present, with a compressive stress in the latitudinal direction. As we move towards the poles only the longitudinal fiber family is present, while at the very poles the principal stresses coincide and the direction of fibers becomes indeterminate.

\subsection{Membrane strip}
\begin{figure}[htp]
\centering
\includegraphics[width=7cm]{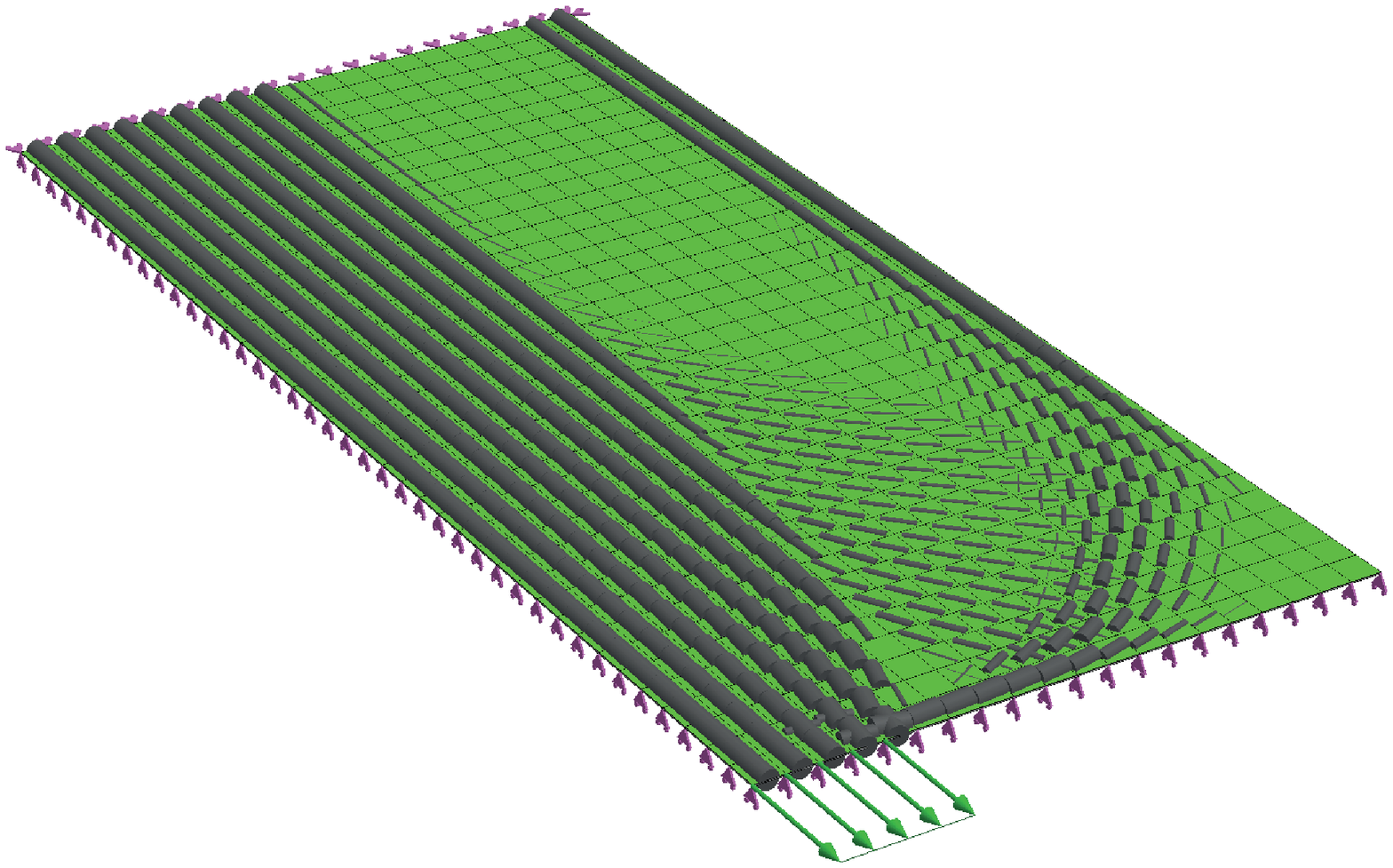}
\includegraphics[width=7cm]{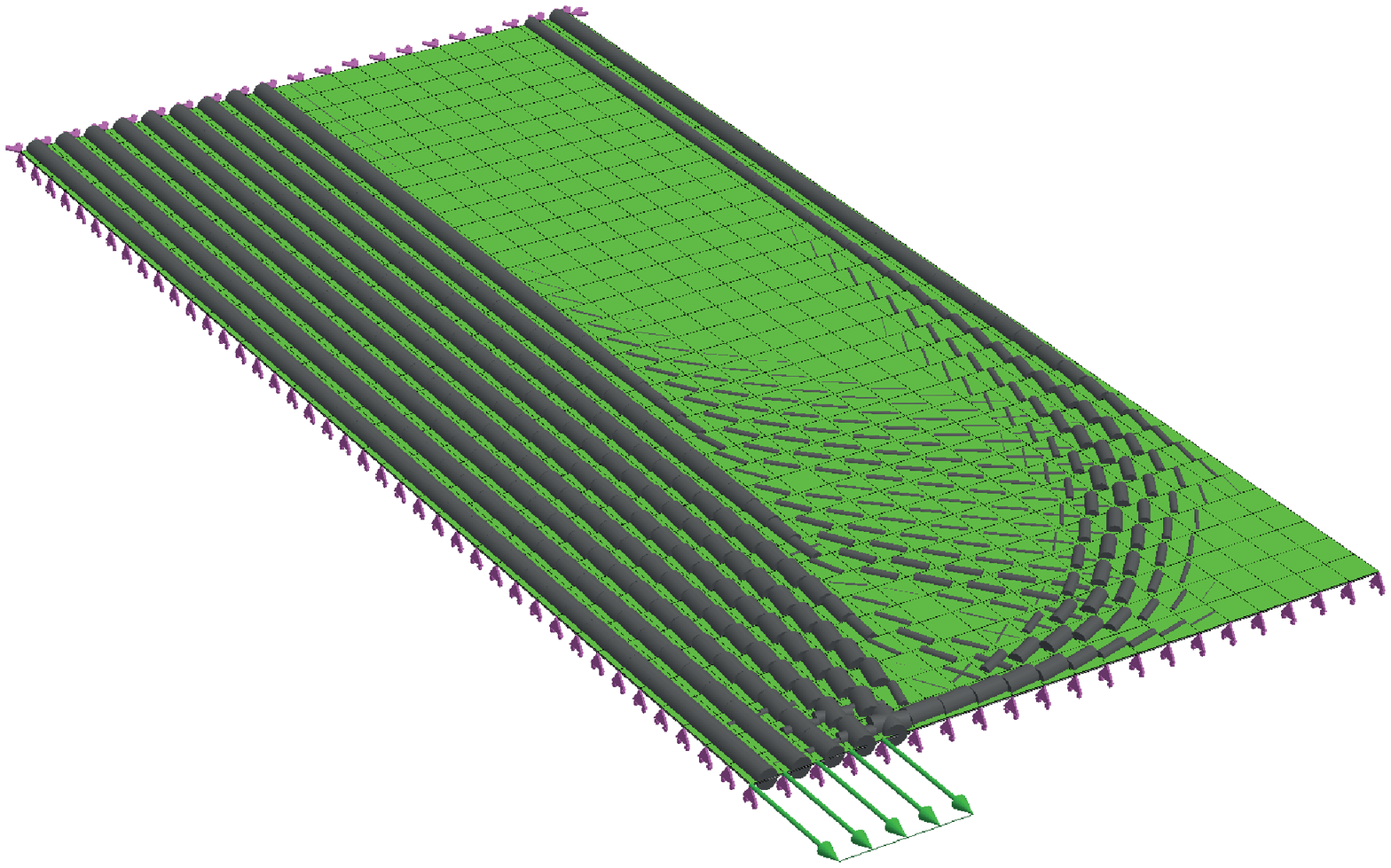}
\caption{Optimal fiber distribution for a rectangular membrane using two different initial fiber directions.}\label{fig:example2}
\end{figure}

A rectangular membrane of shape 1 $\times$ 0.5 is fixed along one of its short sides and loaded by a force $q$ per unit length on a part of length 0.1 of the other short side, as shown in Figure \ref{fig:example2}. The date are $E=1$, $\nu=0$, $t_b=0.005$, $q=0.001$, $V=0.01$, $\overline{t}_1=\overline{t}_2=0.008$, $\underline{t}_1=\underline{t}_2=0$ and $\alpha=2$. As in the previous example, the initial fiber thickness is uniform and chosen so as to satisfy the volume constraint as an equality.

The left hand solution of Figure \ref{fig:example2} is found using initial fiber directions defined by the rectangle sides. The right had solution, on the other hand, uses initial directions defined by principal stress directions found in an initial calculation where fibers are excluded. The left hand problem converged, using the same tolerances as in the previous problem, in 28 optimality criteria iterations and 12 updates of the fiber directions. The right hand problem converged in 15 optimality criteria iterations and 6 updates of the fiber directions. The slightly difference between the two solutions is likely the result of a possible non-uniqueness of the solution of problem $(\mathbb{P})$. However, the objective function values for the two cases are essentially the same.

\section{Conclusions}
The classical facet approach to membrane shells was recently extended to curved elements by Hansbo and Larson \cite{HL}. Here we make a further extension by showing how orthotropic material, of fiber type, can be treated in a similar way, partly inspired by exact plate theory of Nardinocchi and Podio-Guidugli \cite{NPG}. Based on this orthotropic membrane shell theory we formulate a stiffness design problem, where we seek an optimal structure by both rotation and sizing of reinforcing fibers. The two design variables - representing rotation and sizing - naturally splits the formulation into two minimum statements, reminiscent of a Nash game, which suggests a sequential numerical treatment, previously used for similar problems \cite{BS}. This type of formulation also makes clear the distinct character of statically determined problems, which occur for large classes of membrane shells \cite{Ciarlet}. For such problems, the material independent stress state implies that the two minimization statements of $(\mathbb{P})$ become decoupled, and since the sizing of fibers is a convex problem, the full problem $(\mathbb{P})$ essentially inherits this property.

The approach presented in this paper has several intriguing extensions, that would be important for applications such as the design of racing boat sails. Inclusion of pre-stress and wrinkling states related to negative stresses are examples of this. Extension to large deformations, based on the model of Hansbo et al. \cite{HLL}, should also be of clear interest.

\begin{acknowledgements}
This research was supported in part by the Swedish Foundation for Strategic Research Grant No. AM13-0029, the Swedish Research Council Grants Nos. 2011-4992 \& 2013-4708, and the Swedish Research Programme Essence.
\end{acknowledgements}

\section*{Appendix}
\appendix
As a special case of the fiber material defined by $\mathbb{E}^{\mbox{\scriptsize memb}}$, consider the orthotropic material consisting of two orthogonal families of mechanically equal fibers, i.e., $\alpha_1=\alpha_2=\alpha$. We will represent the constitutive law of such a material in the orthogonal base $\{\bm{s},\bm{s}^\bot,\bm{n}\}$, where $\bm{s}=\bm{s}_1$ and $\bm{s}^\bot=\bm{s}_2$. The non-zero part of the stress tensor is $\bm{P}_\Sigma\bm{\sigma}\bm{P}_\Sigma$ and in the indicated base we have:
\begin{equation}\label{eq:stress_11}
  \sigma_{11}:=\bm{s}\cdot(\bm{P}_\Sigma\bm{\sigma}\bm{P}_\Sigma)\bm{s}=\bm{S}:(\bm{P}_\Sigma\bm{\sigma}\bm{P}_\Sigma)=A\varepsilon_{11}+B\varepsilon_{22},
\end{equation}
\begin{equation}\label{eq:stress_22}
  \sigma_{22}:=\bm{s}^\bot\cdot(\bm{P}_\Sigma\bm{\sigma}\bm{P}_\Sigma)\bm{s}^\bot=C\varepsilon_{22}+B\varepsilon_{11},
\end{equation}
\begin{equation}\label{eq:stress_12}
  \sigma_{12}:=\bm{s}\cdot(\bm{P}_\Sigma\bm{\sigma}\bm{P}_\Sigma)\bm{s}^\bot=D(\varepsilon_{12}+\varepsilon_{21})=2D\varepsilon_{12},
\end{equation}
where
\begin{gather*}
\varepsilon_{11}=\bm{s}\cdot(\bm{\varepsilon}\bm{s})=\bm{S}:\bm{\varepsilon},\quad \varepsilon_{22}=\bm{s}^\bot\cdot(\bm{\varepsilon}\bm{s}^\bot),\\ \varepsilon_{12}=\bm{s}\cdot(\bm{\varepsilon}\bm{s}^\bot),\quad
\varepsilon_{21}=\bm{s}^\bot\cdot(\bm{\varepsilon}\bm{s})
\end{gather*}
and
\[
A=\frac{t_b}{t}(\delta+2\mu)+\frac{t_1}{t}\alpha,\quad B=\frac{t_b}{t}\delta,
\]
\[
C=\frac{t_b}{t}(\delta+2\mu)+\frac{t_2}{t}\alpha,\quad D=\frac{t_b}{t}\mu.
\]
Since there is no coupling between normal and shear stresses, one concludes that the principal material directions are given by $\bm{s}$ and $\bm{s}^\bot$. Moreover, the condition defining a so-called low shear material is that the constant $\beta$ below is non-negative, which is indeed the case:
\[
\beta=A+C-2B-4D=\frac{t_1+t_2}{t}\alpha \geq 0.
\]
Moreover, $A \geq C$ obviously follows from $t_1 \geq t_2$.

\end{document}